\documentclass[11pt,letterpaper]{article}

\usepackage[T1]{fontenc}
\usepackage[utf8]{inputenc}
\usepackage{lmodern}
\usepackage[margin=1in]{geometry}
\usepackage{amsmath,amssymb,amsthm}
\usepackage{microtype}
\usepackage[hidelinks]{hyperref}

\numberwithin{equation}{section}
\newtheorem{theorem}{Theorem}[section]
\newtheorem{lemma}[theorem]{Lemma}
\newtheorem{proposition}[theorem]{Proposition}
\newtheorem{corollary}[theorem]{Corollary}
\theoremstyle{definition}
\newtheorem{definition}[theorem]{Definition}
\newtheorem{example}[theorem]{Example}

\newcommand{\Real}{\mathbb{R}}
\DeclareMathOperator{\Ric}{Ric}
\DeclareMathOperator{\Sn}{Sn}
\DeclareMathOperator{\Cn}{Cn}
\newcommand{\abs}[1]{\left\vert#1\right\vert}

\title{Extrinsic Representation of the Partial Ricci Flow}
\author{%
Sergey Liflandsky\thanks{Department of Mathematics and Computer Science,
The Open University of Israel, Haifa Campus, Israel.
\newline E-mail: \href{mailto:sliflandsky@gmail.com}{\texttt{sliflandsky@gmail.com}}.}
\and
Vladimir Rovenski\thanks{Department of Mathematics, University of Haifa,
Mount Carmel, 3498838 Haifa, Israel.
\newline E-mail: \href{mailto:vrovenski@univ.haifa.ac.il}{\texttt{vrovenski@univ.haifa.ac.il}}.}}
\date{}
\hypersetup{pdftitle={Extrinsic representation of the partial Ricci flow},
  pdfauthor={Sergey Liflandsky and Vladimir Rovenski}}

\begin{document}
\maketitle

\begin{abstract}

The partial Ricci flow (PRF) modifies the Riemannian metric on a foliation on a distribution orthogonal to the leaves.
In this work we develop extrinsic representations of PRF on geodesic foliations (twisted/warped products, weak $K$-contact and weak $\beta$-Kenmotsu manifolds) and analyze model cases where PRF admits closed-form solutions. Using the convergence of PRF on model cases, we construct its extrinsic representations initialized by rotational hypersurfaces in Euclidean space as well as parallel curved surfaces and ellipsoids.
The extrinsic viewpoint for PRF provides new tools for studying curvature-driven dynamics on foliations, with potential applications to calculus of moving surfaces
and mathematical physics.

\par\medskip
\noindent\textit{2020 Mathematics Subject Classification:} 53E20 (primary); 35K05, 53C12, 53A05 (secondary).
\par\smallskip
\noindent\textit{Keywords:}
Isometric immersion,
$K$-contact and $\beta$-Kenmotsu structures,
parallel curved surface,
partial Ricci flow,
rotational hypersurface.
\end{abstract}

\tableofcontents

\section{Introduction}

By \textit{geometric flow} we mean the evolution of a metric $g$ on a manifold under a differential equation, $\partial_t g=B(g)$, where the tensor $B$ is associated with curvature.
The~Ricci flow, and the mean curvature flow are most popular geometric flows.
In~\cite{r24new,rov-117,Rov-Wa-2021} geometric flows on a foliated manifold are studied, the metric $g$ is modified along the normal to the leaves distribution,
and the evolution operator is elliptic along the leaves.
The \textit{partial Ricci flow} (PRF), on a foliated manifold $(M,{\cal F})$
\begin{align}\label{E-flow-a}
 \partial_t\,g=-2\Ric^\perp_{g} + 2\,{r}\,g^\perp, \qquad
 {r}\in\mathbb{R}
\end{align}
where $\Ric^\perp_g$ is the partial Ricci curvature,
and $g^\perp$ is the restriction of $g$ on the distribution ${\cal D}$ orthogonal to $\cal F$,
is the tool to prescribe the mixed sectional curvature,
$K_{\rm mix}$. The~PRF preserves ${\cal D}$ and such geometric properties of the leaves as totally umbilical, geodesic or harmonic (minimal) submanifolds.
 The~symmetric (0,2)-tensor
 $\Ric^\bot_g(X,Y)=\sum\nolimits_{\,a} g({R}(X^\bot,\,{e}_a)\,Y^\bot,\,{e}_a)$,
where $\{e_a\}$ is any local orthonormal basis of $\cal D$, is called the \textit{partial Ricci curvature} of
a foliation $\cal F$ of $(M,g)$.
For a one-dimensional foliation $\cal F$
and a unit vector $\xi\in T{\cal F}$, we~have
 $\Ric_g^\perp(X,Y)=g(R_{X,\,\xi}\,\xi,\,Y) =g(R_{\xi}\,X,\,Y)$,
where $R_\xi$ is the \textit{Jacobi operator} for~$\xi$.
For a foliated surface $(M^2,{\cal F},g)$,
we have $\Ric^\perp_{g}=K\,g^\perp$ and $R_\xi=K\,I_{\cal D}$, where $K$ is the gaussian curvature,
hence, \eqref{E-flow-a} reduces to the 
$\cal D$-conformal flow
 $\partial_t\,g = -2\,(K-{r})\,g^\perp$.

For better understanding and visualization the geometric flows, we can represent them extrinsically by depicting a submanifold deforming in Euclidean space,
when the metric of the flow
is the metric induced by the manifold's~embedding.
In this paper, we adapt to PRF the following concept, studied in \cite{CDJ-2020} for the Ricci flow.

An \textit{extrinsic representation} ensures that, given a geometric flow $g_t$ on a differentiable manifold $M$, there is a smooth family of immersions $\psi_t:M \to S_t\ (t\in\mathbb{R})$ in a fixed Euclidean space, such that the metric extrinsically induced on $S_t$ by the ambient Euclidean metric is precisely~$g_t$.

This approach is motivated by Nash's embedding theorems, which state that any Riemannian manifold can be isometrically embedded into some Euclidean space.
Moreover, for a short time, one can construct a smooth one-parameter family of isometric embeddings of any smooth parametrized family of Riemannian metrics on a smooth closed manifold into Euclidean space, see~\cite{Zeg-17}, in particular, the Ricci flow on
$M^n$ can be represented in $\mathbb{R}^m$, where
$m=\max\{\frac12\,n(n+5),\, \frac12\,(n+3)+5\}$.

The PRF exists for a short time, see \cite{Rov-Wa-2021}, so \textit{what is the relation between $m$ and $n$ for the extrinsic representations of the {\rm PRF} on a smooth foliated manifold $M^n$ in~$\mathbb{R}^m$}?

Some~authors studied extrinsic representations of the Ricci flow initiated by rotational hypersurfaces in space forms, e.g., \cite{CDJ-2020,Suda-2024},
whose metrics are rotational symmetric -- a special class of warped products.
Meridian curves of a rotational hypersurface form a geodesic foliation.
Such foliations arise naturally on twisted and warped products $B\times_{\rho}\bar M$, on weak $\beta$-Kenmotsu manifolds and {weak contact metric manifolds} $M(f,\xi,\eta,g)$, see \cite{RWo-2,rov-126,rov-117},
obtained by replacing the complex structure of the distribution $\ker\eta$ with a non-singular skew-symmetric tensor.
The~main achievement of the PRF for geodesic foliations
is the dynamical characterizations of $K$-contact
manifolds and warped products of constant mixed curvature, $K_{\rm mix}$

--
PRF evolves a weak $K$-contact structure (with $K_{\rm mix}>0)$, exponentially fast towards a $K$-contact structure (with $K_{\rm mix}\equiv1)$

--
PRF evolves warped products toward warped products with
$K_{\rm mix}={\rm const}$.
This has practical implications because a simply connected weak Kenmotsu $\beta$-manifold with a complete Reeb field $\xi$ such that $\nabla\beta$ is collinear with $\xi$ is a warped product.

These advances raise the following problem:
\textit{to study extrinsic representations of the {\rm PRF} on geodesic foliations initialized by rotational hypersurfaces in Euclidean space, or hypersurfaces from other distinguished classes}.

The paper focuses on extrinsic representations of PRF
on geodesic foliations, when $g_t$ correspond to
a warped product, $\beta$-Kenmotsu or $K$-contact structures, and $S_t$ are rotational hypersurfaces or ellipsoids, respectively.
We distinguish 3 assertions: existence of the PRF, realization by a certain isometric immersion, and smooth extension of such immersion across a pole. The first doesn't imply the other two.

In Section~\ref{sec:02} we show
the convergence of the PRF for some warped and twisted products, and give examples.
In Section~\ref{sec:03} we study the PRF initiated by rotational hypersurfaces in Euclidean space
and by parallel curved surfaces.
In Section~\ref{sec:04} we recall the notion of
weak almost contact metric manifolds and study extrinsic representations of the PRF
on weak $K$-contact manifolds initiated by ellipsoids,
and on weak $\beta$-Kenmotsu manifolds initiated by rotational hypersurfaces.

\section{The PRF on Twisted and Warped Products}
\label{sec:02}

The notion of a warped and twisted product is popular in  differential geometry as well as in general relativity (some spacetime models are warped products).

Let $(B^s, g_B)$ be a compact Riemannian manifold with boundary, and $(\bar M, \bar g)$ be an $n$-dimensional Riemannian manifold.
The~\textit{twisted product} $B\times_{\rho}\bar M$
is the manifold $M=B\times\bar M$ with the metric
$g=g_B+\rho^2\,\bar g$,
where $\rho\in C^\infty(M)$ is a positive function.
When $\rho$ depends on $B$ only, the twisted product is a~\textit{warped product}.
The~special case of a warped product are \textit{rotational symmetric metrics}, i.e., $(\bar M,\bar g)$ is a unit $n$-sphere, such metrics appear on rotational hypersurfaces in space~forms.

The \textit{leaves} $B\times\{v\}$ of a twisted product
are totally geodesic, and the \textit{fibers} $\{u\}\times\bar M$ are totally umbilical with the second fundamental form
$h=-\nabla^{\cal F}(\log\rho)\,\bar g$.
For twisted products $B^s\times_{\rho}\bar M^n$, we have
$\Ric^\perp_{g}=(S^g_{\rm mix}/n)\,g^\perp$,
where $S^g_{\rm mix}=\operatorname{tr}\Ric^\perp_{g}$,
and the PRF becomes a ${{\cal D}}$-conformal flow
 $\partial_t\,g=-2\,(S^g_{\rm mix}-{r})\,g^\perp$.

Since
$R(X,\xi,Y,\xi)=-\rho^{-1} \xi(\xi(\rho))\,\bar g(X,Y)$ for unit $\xi\in T{\cal F}$,
the mixed sectional curvature $K_{\rm mix}$ depends on a point of $B$ only, and
 $\Ric^\perp_g=-(\Delta_{\cal F}\,\rho/\rho)g^\perp$.
Here, $\nabla^{\cal F}$ and $\Delta_{\cal F}$ denote,
respectively, the gradient and the laplacian on the leaves.

If twisted product metrics $g_t=g_B+\rho^2(t,u,v)\bar g$
solve the PRF on $B\times\bar M$, then
 $\partial_t(\rho^2)\,\bar g = 2\,(\Delta_{\cal F}\,\rho/\rho)\,\rho^2\bar g
 +2\,{r}\,\rho^2\,\bar g$.
This yields the leaf-wise parabolic PDE
\begin{equation*}
\partial_t\,\rho=\Delta_{\cal F}\,\rho+{r}\,\rho.
\end{equation*}
Let $\mu(t,u)=\rho(t,u)\,|_{\,\partial B}$ be twice continuously differentiable in $t$, and
\begin{equation}\label{E-mu-twice}
\lim\limits_{\,t\rightarrow\infty}\mu(t,u)=\tilde\mu(u),\qquad
 \lim\limits_{\,t\rightarrow\infty}\partial_t\mu(t,u)=0,\qquad
 \lim\limits_{\,t\rightarrow\infty}\partial^2_t\mu(t,\,\cdot)=0
\end{equation}
for a smooth nonnegative function $\tilde\mu:\partial B\to\Real$ uniformly for $u\in\partial B$.

Define
\begin{align*}
\delta(t,\,\cdot)&=\mu(t,\,\cdot)-\tilde\mu,\\
\nu(t)&=\max
 \big\{\,\|\,\delta(t,\,\cdot)\,\|_{\,C^0(\partial B)},\;\|\,\partial_t\delta(t,\,\cdot)\,\|_{\,C^0(\partial B)}\,\big\},\\
\tilde\nu(t)&=\max
 \big\{\,\|\,\delta(t,\,\cdot)\,\|_{\,C^0(\partial B)},\;\|\,\partial_t\delta(t,\,\cdot)\,\|_{\,C^0(\partial
 B)}, \;\|\,\partial^2_{t,t}\delta(t,\,\cdot)\,\|_{\,C^0(\partial B)}\,\big\}.
\end{align*}

We first record the fact identifying the PRF with the profile heat equation, it is the hinge on which the extrinsic representation turns.

\begin{theorem}[Rovenski-Zelenko]\label{T-stabnonstat}
Let the warped product metrics $g_t=g_B+\rho^2(t,u) \bar g$ on $B^p\times_{\rho}\bar M^n$ solve the {\rm PRF} \eqref{E-flow-a}, and any of conditions {\rm i)}--{\rm iii)} are satisfied

\quad {\rm i)}
${r}<0$ and \eqref{E-mu-twice}$_{1,2}$\qquad
{\rm ii)} $0\le{r}<\lambda_1$, $p<4$, and \eqref{E-mu-twice}

\quad {\rm iii)} ${r}=\lambda_1$, $p<4$, \eqref{E-mu-twice} and $\tilde\mu=0$,
 \ $\int_0^\infty\nu(\tau)\,{\rm d}\tau<\infty$.

Then $g_t$ exist for all $t\ge0$, and converge, as $t\to\infty$, uniformly on $B^p\times\bar M^n$ in $C^0$-norm
to the limit warped product metric $g_\infty=
g_B+\rho_\infty^2(u)\,\bar g$ with
$\Ric_{g_\infty}^\perp=r\,\bar g$.
Moreover, in cases {\rm i)} and {\rm ii)}, the {\rm PRF} has a global single point attractor, while in case {\rm iii)} the limit metric $g_\infty$ depends on initial and boundary conditions.
\end{theorem}

\begin{proof}
Let $0<\lambda_1\le\ldots\le\lambda_i\le\ldots$ be eigenvalues of the Schr\"odinger operator
 $-\Delta_{\cal F}-{r}$
on~the leaves $B\times\{y\}$ counting their multiplicities, $\lim_{\,j\to\infty}\lambda_j=\infty$,
and $\{e_i\}_{1\le i<\infty}$ corresponding $L_2$-orthonormal basis of eigenfunctions
with $e_{i}=0$ on $\partial B$ of the eigenvalue problem:
 $-\Delta_{\cal F}\,e_i=\lambda_i\,e_i$.
Then, $\lambda_1$ has multiplicity 1 (hence, $\lambda_1<\lambda_2$) and $e_1>0$ on the interior of $B$.
The leaf-wise parabolic Cauchy's problem
\begin{equation}\label{E-r}
 \partial_t\,\rho=\Delta_{\cal F}\,\rho+{r}\,\rho,\qquad
 \rho(0,\,\cdot)=\rho_0,\qquad
 \rho(t,\,\cdot)\,|_{\,\partial B}=\mu(t,\,\cdot)
\end{equation}
has a~unique solution $\rho:[0,\infty)\times B\to\Real$.
One can prove convergence of $\rho$, as $t\to\infty$, to a~stationary state -- a~solution
$\rho_\infty : B\to\Real$ of the problem
\begin{equation}\label{statprob}
 -\Delta_{\cal F}\,\rho_\infty = {r}\,\rho_\infty,\qquad
 \rho_\infty\,|_{\,\partial B}=\tilde\mu.
\end{equation}
The rest of proof is similar to the proof of
\cite[Theorem~5.98]{Rov-Wa-2021}.
\end{proof}

For $p=1$ and $B=[0,1]$, the problems \eqref{E-r} and \eqref{statprob} have the following view
\begin{align}\label{eq:prfheat}
 \rho'_{t}=\rho''_{u,u}+r\,\rho ,\qquad
 \rho(0,\,\cdot)=\rho_0,\qquad
 \rho(t,\,0)=\mu_0(t),\qquad \rho(t,\,1)=\mu_1(t) \\
\label{eq:statbvp}
 \rho''_\infty+r\,\rho_\infty=0\ \ \text{on}\ (0,1),\qquad
\rho_\infty(0)=\alpha_1\in\Real,\qquad \rho_\infty(1)=\alpha_2\in\Real.
\end{align}

To describe stationary profiles for $r\in\Real$, we use the following pair of solutions of $y''+r\,y=0$ normalized by $y(0)=0$, $y'(0)=1$ and $y(0)=1$, $y'(0)=0$, respectively,

\begin{equation*}
\Sn_{r}(u)=\begin{cases}
{\sin(\sqrt{r}\,u)}/{\sqrt{r}} &r>0\\[2pt] u & r=0 \\[2pt] {\sinh(\sqrt{-r}\,u)}/{\sqrt{-r}} &r<0 \end{cases}\qquad
\Cn_{r}(u)
=\begin{cases}\cos(\sqrt{r}\,u) &r>0\\[2pt] 1 & r=0 \\[2pt] \cosh(\sqrt{-r}\,u) &r<0.
\end{cases}
\end{equation*}

The functions $\Sn_{r}$, $\Cn_{r}$ are jointly smooth in $(u,r)$ and satisfy the following relations:
$\Sn_{r}'=\Cn_{r},\ \Cn_{r}'=-r\,\Sn_{r},\ \Cn_{r}^{2}+r\Sn_{r}^{2}=1$.

\begin{corollary}[Stationary profiles]\label{prop:stationary}
Let $r<\pi^{2}$ and $\alpha_1,\alpha_2\in\Real$.
The problem \eqref{eq:statbvp} has a unique solution
$\rho_\infty(u,r;\alpha_1,\alpha_2)
 =\dfrac{\alpha_2\Sn_{r}(u)+\alpha_1\Sn_{r}(1-u)}{\Sn_{r}(1)}$
 linear\ in $\alpha_1,\ \alpha_2$.
For $r=\pi^{2}$ and $\alpha_1,\alpha_2\ge0$ with $\alpha_1^{2}+\alpha_2^{2}>0$, \eqref{eq:statbvp} is solvable if and only if $\alpha_2=-\alpha_1$, in which case solutions are $\rho_\infty(u,r;\alpha_1,\alpha_2)=C\sin(\pi u)$ for $C\in\Real$.
\end{corollary}

\begin{example}\label{lem:prf-heat}
i) For $M=[0,1]\times{\mathbb S^{1}}$ 
the twisted product (semi-geodesic) metrics are
\begin{equation}\label{eq:semigeo}
g=\,{\rm d}u^{2}+\rho^{2}(u,v)\,{\rm d}v^{2},\qquad \rho>0 .
\end{equation}
For the metric \eqref{eq:semigeo},
each $u$-curve (i.e., $v=\mathrm{const}$) is a unit-speed geodesic, and the Gaussian curvature
is $K=-{\rho''_{u,u}}/{\rho}$ (involves no $v$-derivatives of $\rho$).

ii) Let $g_{t}=\,{\rm d}u^{2}+\rho^{2}(t,u) \,g_{\mathbb{S}^{n-1}}$,
$\rho(t,\cdot)>0$ on $(0,1)$, be rotationally symmetric metrics on $[0,1]\times\mathbb{S}^{n-1}$. Then
$g_{t}$ solve the PRF \eqref{E-flow-a} if and only~if
\eqref{eq:prfheat}~holds.

Indeed, the leaves $\{\omega={\rm const}\}$ compose a geodesic foliation ${\cal F}$ with unit tangent $N=\partial_u$ and leaf-wise Laplacian
$\Delta_{\cal F}=\partial_{u,u}^{2}$,
and ${\cal D}=T\,\mathbb{S}^{n-1}$,
this is the warped product $[0,1]\times_{\rho}\mathbb{S}^{n-1}$.
We obtain the formula
 $R_N(t) =-(\rho''_{u,u}/\rho)\,g_t^{\perp}$,
while $g_t^{\perp}=\rho^{2}g_{\mathbb{S}^{n-1}}$
and $\partial_t g_t=(\rho^{2})'_{t}\,g_{\mathbb{S}^{n-1}}=2\,\rho\,\rho'_{t}\,g_{\mathbb{S}^{n-1}}$.
Thus, the PRF \eqref{E-flow-a}~reads
 $2\,\rho\,\rho'_{t}=2\,(\rho''_{u,u}/\rho)\rho^{2}+2\,r\,\rho^{2}$,
that reduces to \eqref{eq:prfheat}.
\end{example}

\begin{proposition}\label{thm:A}
Consider metric $g_{0}={\rm d}u^{2} +\rho_{0}^2(u,v)\,{\rm d}v^{2}$ on $[0,1]\times{\mathbb S}^1$ with smooth $\rho_{0}>0$ on $(0,1)\times{\mathbb S}^1$, with boundary values $\alpha_1(v):=\rho_{0}(0,v)\ge0$ and $\alpha_2(v):=\rho_{0}(1,v)\ge0$ fixed in time. Then the following~holds

{\rm i)} the {\rm PRF} with the boundary conditions
 $\rho(t,0,v)=\alpha_1(v)$
 and
  $\rho(t,1,v)=\alpha_2(v)$
has a unique solution $g_{t}=\,{\rm d}u^{2}+\rho^{2}(t,u,v)\,{\rm d}v^{2}$ for all $t\ge0$, it is given for $r<\lambda_{1}$~by

\begin{equation*}
\rho(t,u,v)=\rho_\infty\big(u,r,\alpha_1(v),\alpha_2(v)\big)
+\sum\nolimits_{\,n=1}^{\infty}{A}_{n}(v)\,\sin\big({n\pi}\,u\big)\exp((r-\lambda_{n})t)
\end{equation*}

where
${A}_{n}(v)={2}\int_{0}^{1}\Big(\rho_{0}(u,v)
-\rho_\infty\big(u,r,\alpha_1(v),\alpha_2(v)\big)\Big) \sin\big({n\pi}\,u\big)\,{\rm d}u$.

The function $\rho>0$ on $(0,1)\times{\mathbb S}^1\times[0,\infty)$, so $g_{t}$ is a Riemannian metric on the interior for all $t$, degenerating at a boundary circle when the boundary value vanishes.

{\rm ii)} if $r<\lambda_{1}$ then,
$\|\,\rho(t,\cdot,v)-\rho_\infty\big(\cdot,r,\alpha_1(v),\alpha_2(v)\big)\|_{\,C^{k}([0,1])}
\le C_{k}\exp((r-\lambda_{1})t)\ (t\ge1)$,
uniformly in $v$ for every $k\ge0$, consequently, $g_{t}\to g_{\infty}$ in every $C^{k}$ norm, where
$g_{\infty}=\,{\rm d}u^{2} +\rho_\infty^{2}\big(u,r,\alpha_1(v),\alpha_2(v)\big)\,{\rm d}v^{2}$
has constant Gaussian curvature $K\equiv r$.

{\rm iii)} if $r=\lambda_{1}$, a limit exists only if $\alpha_1\equiv\alpha_2\equiv0$, in that case

\begin{equation*}
\lim_{t\to\infty}\rho(t,u,v)=B_{1}(v)
\sin\big({\pi}u\big),\qquad
B_{1}(v)={2}\int_{0}^{1}\rho_{0}(u,v) \sin\big({\pi}u\big)\,{\rm d}u>0
\end{equation*}

in any $C^{k}$ norm, and the limit metric has constant curvature $r$ on the interior.

{\rm iv)} if $r>\lambda_{1}$, then $\lim_{\,t\to\infty}g_{t}$ exists if and only if for any $v$ and $n$ with $\lambda_{n}\le r$, the coefficient ${A}_{n}(v)$ of the solution vanishes (with $\rho_\infty$ replaced by any particular stationary solution when one exists), for generic data no limit~exists.
\end{proposition}

\begin{proof} Fourier series are used here in the standard way.
\end{proof}

\section{The PRF on Rotational and Parallel Curved Surfaces}
\label{sec:03}

In this section, we construct an extrinsic representation of the PRF \eqref{E-flow-a} initialized by rotational hypersurfaces and parallel curved surfaces of Euclidean space.
\subsection{Rotational Hypersurfaces}

Rotational hypersurfaces with a one-dimensional axis in
$\mathbb{R}^{n+1}$ are foliated into $(n-1)$-spheres~(parallels) and equipped with rotationally symmetric metrics  $g={\rm d}u^2+\rho^2(u)\,{\rm d}s^2_{n-1}$ -- a~special case of a warped product.
Such a hypersurface (and~its profile) can be given as a graph $\rho=F(x_{n+1})$, where $F\in C^2$,
$\,\rho=\sqrt{x_1^2+\ldots+x_{n}^2}$
and $(x_i)$ are Cartesian coordinates in~$\mathbb{R}^{n+1}$.
We assume $F''\not\equiv0$ to exclude the cases
when the profile curve is a straight,
and $M$ is a cone, a~cylinder, or a hyperplane.

A rotational hypersurface in Euclidean space $\mathbb{R}^{n+1}$ with a multi-dimensional axis is a~hypersurface obtained by rotating a lower-dimensional profile around a linear subspace (the ``axis'') of dimension $>1$ (not a line).
So, instead of the rotational hypersurface with axis = line, we get a foliated hypersurface
invariant under the group
$O(n+1-s)$ acting on the orthogonal complement of the $s$-dimensional axis: for each $y\in\mathbb{R}^s$, the fiber in $x$ is a round $(n-s)$-sphere of radius $\rho=F(y)$ in~$\mathbb{R}^{n+1-s}$.

\begin{proposition}\label{prop:rbound}
Let $S^n\subset\mathbb{R}^{n+1}$ be a
rotational hypersurface with a profile curve $(\rho,h)$ parameterized by arclength, where $\rho,h:[0,1]\to\Real$ are smooth and $\rho>0$ on $(0,1)$.
If $S^n$ has $K_{\rm mix}=r={\rm const}>0$,
then $r\le\pi^{2}$, and for $r=\pi^{2}$, we have $\rho(0)=\rho(1)=0$, and for any $C\in(0,{1}/{\pi}]$, the functions $(\rho,h)$ are given~by
\begin{align}\label{eq:piL}
\rho(u)=C\sin\big(\pi\,u\big),
\qquad
h(u)=\int_{0}^{u}\sqrt{1-\big(\pi\,C\big)^{2} \cos^{2}\big(\pi\,s\big)}\,{\rm d}s .
\end{align}
\end{proposition}

\begin{proof} This is standard.
\end{proof}

\begin{example}[see Example~5.101 in \cite{Rov-Wa-2021}]
\label{Ex-tilder}\rm
 Consider rotational surfaces $S_t\subset\Real^3$
\[
 {\bf r}=\big(\rho(t,u)\cos v,\ \rho(t,u)\sin v,\ h(t,u)\big),\qquad
 0\le u\le 1,\qquad |v|\le\pi,\qquad \rho\ge0
\]
with $(\rho'_{u})^2+(h'_{u})^2=1$.
Profile curves $v={\rm const}$ of $S_t$ are unit speed geodesics.
The~$v$-curves of $S_t$ are circles (in~$\Real^3$) of geodesic curvature $k=-(\log\rho)'_{\,u}$.
The~metric $g_t ={\rm d}u^2+\rho^2(t,u)\,{\rm d}v^2$ is {rotational symmetric}, and its Gaussian curvature is
$K=-\rho''_{\,u,u}/\rho$.
Let $g_t$ solve the PRF, then $\rho(t,u)$ solves the
problem \eqref{eq:prfheat}.

Assume that $\mu_0,\mu_1\in C^1[0,\infty)$,
then there exist limits
 $\lim\limits_{\,t\to\infty}\mu_j(t)=\tilde\mu_j\in [0,\infty)$.

The solution of stationary problem \eqref{eq:statbvp}
with $B=[0,1]$ and $\lambda_j=(\pi\,j)^2$ is
\begin{equation}\label{E-stat-rho}
 \rho_\infty(u)=\left\{\begin{array}{cc}
 \dfrac{\tilde\mu_1\sin(\sqrt{r}\,u) +\,\tilde\mu_0\sin(\sqrt{r}\,(1-u))}{\sin(\sqrt{r})} &
 {\rm if} \ \ 0<{r}<\lambda_1\\
 \tilde\mu_0+(\tilde\mu_1-\tilde\mu_0)(u) & \quad {\rm if}\ {r}=0\\
 \dfrac{\tilde\mu_1\sinh(\sqrt{-{r}}\,u) +\,\tilde\mu_0\sinh(\sqrt{-{r}}\,(1-u))}{\sinh(\sqrt{-{r}})} &
 \quad {\rm if}\ {r}<0.\\
 \end{array}\right.
\end{equation}
For the resonance case, ${r}=\lambda_1=\pi^2$, the stationary problem is solvable if and only if
$\tilde\mu_0=\tilde\mu_1=0$,
and the solution is $\rho_\infty(u)=C\sin(\pi u)$, where $C={\rm const}>0$.

If~${r}>\pi^2$, then $g_t$ diverge as $t\to\infty$,
and if~${r}\le\pi^2$, then $g_t$ converge to the limit metric $g_\infty={\rm d}u^2+\rho\,_\infty^2(u)\,{\rm d}v^2$
with $K(g_\infty)={r}$, see Theorem~\ref{T-stabnonstat}.

With~known solution $\rho(t,u)\ (t\ge0)$ of \eqref{eq:prfheat}
and $|\,\rho'_{u}|\le1$, we get
$h=h(t,0)+\int_0^u\!\sqrt{1-(\rho'_{u})^2}\,{\rm d}u$,
that correspond to rotational surfaces in $\Real^3$ of positive constant Gaussian curvature.
 Assume for simplicity that $\mu_j(t)\equiv\tilde\mu_j$. Then

\begin{equation*}
\tilde\mu_j=0,\qquad
 \int_0^\infty|\mu_j(\tau)-\tilde\mu_j|\,{\rm d}\tau<\infty,\qquad
 \int_0^\infty|\mu^\prime_j(\tau)|\,{\rm d}\tau<\infty, \qquad
 j=0,1.
\end{equation*}

For ${r}=\pi^2$, we get $\rho_\infty=C\sin(\pi\,u)$, where $C=v_1^0=\int_0^1 v_0(s)\,{\rm d}s$,
and then
\begin{equation*}
h_\infty=\frac 1\pi\,{\rm EllipticE}\,(\cos(\pi u), C\pi).
\end{equation*}
Here,
\begin{equation*}
{\rm EllipticE}(z,k)=\int_0^z \sqrt{(1-k^2 s^2)/(1-s^2)}\,{\rm d}s
\end{equation*}
is the incomplete elliptic integral.
To~provide a numerical example for $r=\pi^2$, let $u\in[0,\pi]$, \ $C=1$ and $\mu_j=0$. Then, the limit profile curve is a semi-circle $[\sin\,u,\cos\,u]$ for $0\le u\le\pi$, and the limit rotational surface in $\Real^3$ is a round~sphere.
\end{example}

Next, we prove the extrinsic realization theorem:
the rotational hypersurfaces built from the normalized solution realize the PRF extrinsically for all $t\ge0$,
and converge, as $t\to\infty$, to the rotational hypersurface of
$K_{\rm mix}={r}$ when ${r}\le\lambda_1$.

The key analytic fact is that the \textit{total slope}
$\delta:=\sup\nolimits_{\,(t,u)\in[0,1]\times[0,\infty)}\abs{\,\rho'_{u}(t,u)}$
is finite whenever $r\le\lambda_{1}$, and one can find an explicit bound of $\delta$.

\begin{definition}\label{defn:normalized}\rm
Let $r\le\lambda_{1}$, and $\rho_0:[0,1]\to\Real$ be smooth and positive on $(0,1)$ with $\rho_0(0)=\mu_0\ge0$, $\rho_0(1)=\mu_1\ge0$ (and $\mu_0=\mu_1=0$ if $r=\lambda_{1}$), and $\delta<\infty$ be the total slope of the solution $\rho$ of \eqref{eq:statbvp}.
We assume $\delta>0$: the degenerate case $\delta=0$ occurs exactly when $\rho\equiv\mu_0=\mu_1>0$ is constant, and then
we get a trivial representation: every $M_t$ below is a fixed round cylinder. The~\emph{normalized {\rm PRF} solution} is $\widehat{\rho}:=\rho/\delta$,
which by linearity of \eqref{eq:statbvp} is the solution with data $\rho_0/\delta$, $\mu_0/\delta$, $\mu_1/\delta$, and satisfies
\begin{equation}\label{eq:normslope}
\sup\nolimits_{\,(t,u)}\abs{\,\widehat{\rho}\,'_{u}(t,u)\,}=1 .
\end{equation}
\end{definition}

We represent the PRF by rotational hypersurfaces in $\Real^{n+1}$.
Given the normalized PRF solution $\widehat{\rho}(t,u)$, define the map $\psi_t:[0,1]\times\mathbb{S}^{n-1}\to\Real^{n+1}$ for $t\ge0$
\begin{align}\label{eq:embedding}
\psi_t(u,\omega)=\big(\,\widehat{\rho}\,(t,u)\,\omega,\ h(t,u)\big),\qquad
h(t,u)=\int_{0}^{u}\sqrt{1-(\,\widehat{\rho}\,'_{s})^2(t,s)}\,{\rm d}s
\end{align}
where $\omega\in\mathbb{S}^{n-1}(1)$,
the square root in the formula for $h(t,u)$ is real due to \eqref{eq:normslope}, and the choice $h(t,0)=0$, $h'_u\ge0$ pins down the vertical translation and reflection
(an extrinsic representation is defined up to rigid motions of $\Real^{n+1}$).

\begin{lemma}\label{lem:embedding}
For each $t\ge0$ the map \eqref{eq:embedding} parameterizes an embedded rotational hypersurface $S_{t}\subset\Real^{n+1}$ {\rm (}immersed open cylinder $(0,1)\times\mathbb{S}^{n-1}$, closing up at a pole $u\in\{0,1\}$, where a boundary radius $\widehat{\rho}$ vanishes, the closing is smooth when $\abs{\,\widehat\rho\,'_{u}}=1$ there, and is conical otherwise{\rm )}, with the following properties

{\rm i)} The induced metric is $g_{t}=\,{\rm d}u^{2} +\widehat{\rho}\,^2(t,u)\,g_{\mathbb{S}^{n-1}}$ and satisfies the {\rm PRF} \eqref{E-flow-a}

{\rm ii)} For each $\omega\in\mathbb{S}^{n-1}$,
$\psi_t(u,\omega)$ in \eqref{eq:embedding} is a unit-speed geodesic of $S_{t}$ of length~$1$.

Moreover, $\psi_t$ is a $C^{1}$-regular immersion for every $t$, smooth away from the closed set $\abs{\,\widehat\rho\,'_u}=1$, the metric $g_t$ is smooth in $(t,u)$ for $t>0$ and continuous at $t=0$.
\end{lemma}

\begin{proof} This is standard.
\end{proof}

Let $S_{t}\ (t\ge0)$ be rotational hypersurfaces with a common axis and profile curves $\gamma_{t}$ parameterized by arclength on the interval $[0,1]$.
Then $S_{t}$ converge to the \emph{limiting hypersurface} with a profile curve $\gamma_{\infty}$, if and only if $\gamma_{t}\to\gamma_{\infty}$ uniformly on~$[0,1]$.

\begin{theorem}[Extrinsic realization]
Let $r\le\lambda_{1}$, with $\mu_0=\mu_1=0$ if $r=\lambda_{1}$, $\rho_0$ be smooth and positive on $(0,1)$ with $\rho_0(0)=\mu_0\ge0$, $\rho_0(1)=\mu_1\ge0$,  $\delta\in(0,\infty)$ be the total slope.
Suppose that $S_{t}\subset\Real^{n+1}$ are rotational hypersurfaces \eqref{eq:embedding} built from the normalized solution $\widehat{\rho}=\rho/\delta$. Then $\{S_{t}\}_{t\ge0}$ realize the
{\rm PRF} \eqref{E-flow-a} extrinsically\textup{:}  the metric induced on $S_t$ is $g_t=\,{\rm d}u^2+\widehat\rho^{\,2}(t,u)\,g_{\mathbb{S}^{n-1}}$, and

{\rm i)} if $r<\lambda_{1}$, the limiting, as $t\to\infty$, 
hypersurface has
$K_{\rm mix}\equiv r$
and its profile is
\begin{equation}\label{eq:limitgen}
\gamma_{\infty}(u)=\Big(\rho_\infty\big(u,r,{\mu_0}/{\delta}, {\mu_1}/{\delta}\big),\
\int_{0}^{u}
\big({1-\rho'_\infty(s,r;{\mu_0}/{\delta}, {\mu_1}/{\delta})^{2}}\big)^{1/2}\,{\rm d}s\Big)
\end{equation}

{\rm ii)} if $r=\lambda_{1}$ (so $\mu_0=\mu_1=0$), the limiting, as $t\to\infty$, hypersurface exists and is given by \eqref{eq:piL} with
$B={A_{1}}/{\delta}\in(0,{1}/{\pi}]$, where
$A_{1}={2}\int_{0}^{1}\rho_0(u)
\sin\big({\pi}\,u\big)\,{\rm d}u>0$

{\rm iii)} if $r>\lambda_{1}$, then for generic data, $\delta=\infty$,
the normalization is impossible, and no limiting hypersurface exists. Indeed, the intrinsic limit fails by Proposition~\ref{thm:A}iv) unless all coefficients $\widetilde{A}_{i}$ with $\lambda_{i}\le r$ vanish.
\end{theorem}

\begin{proof}
The induced metrics are
$g_t=\,{\rm d}u^2+\widehat\rho^{\,2}g_{\mathbb{S}^{n-1}}$,
they solve the PRF \eqref{E-flow-a}, and the meridians are geodesics of length $1$ for every $t$, see Lemma~\ref{lem:embedding}\,i),ii) and Example~\ref{lem:prf-heat}\,ii).
Note that $\,\widehat{\rho}$ is the solution of \eqref{eq:statbvp} with data $(\rho_0/\delta, \mu_0/\delta, \mu_1/\delta)$ by linearity, unique and smooth for $t>0$ by \cite[Section~5.4.4]{Rov-Wa-2021}
and \eqref{eq:normslope} makes \eqref{eq:embedding} well defined for all $t$. Thus, $\psi_t$ are the required immersions, smooth for $t>0$.

i) By Proposition~\ref{thm:A}\,ii) with the normalized data, $\widehat{\rho}(t,\cdot)\to\rho_\infty(\cdot\,,r,\mu_0/\delta,\mu_1/\delta)$ in $C^{1}([0,1])$, as $t\to\infty$, using
$\rho_\infty(\cdot\,,\mu_0/\delta,\mu_1/\delta)
=\tfrac{1}{\delta}\,\rho_\infty(\cdot\,,\mu_0,\mu_1)$.
Passing to the limit in $\abs{\,\widehat{\rho}\,'_{u}}\le1$ gives  
the derivative of $\rho_\infty$, see \eqref{eq:slopecond} and also Corollary~\ref{prop:stationary}.~Then
\begin{equation*}
 \abs{h(t,u)-h_{\infty}(u)}
 \le\int_{0}^{1}\big|\sqrt{1-(\,\widehat{\rho}\,'_{s})^2(t,s)}
 -\sqrt{1-(\,\rho'_{\infty,u})^2(s)}\,\big|\,{\rm d}s \rightarrow0
\end{equation*}
uniformly in $u$, since
$\widehat{\rho}\,'_{u}\to\rho'_\infty$ uniformly with values in $[-1,1]$. Hence, the profile curves $(\,\widehat{\rho},h)(t,\cdot)$ converge uniformly to \eqref{eq:limitgen}, which is the limiting hypersurface.
A~limiting metric is $g_\infty=\,{\rm d}u^2 +\rho_\infty^2(u)\,
g_{\mathbb{S}^{n-1}}$, and any plane containing $\partial_u$ has sectional curvature
$K_{\rm mix}=-{\rho''_\infty}/{\rho_\infty}=r$
by \eqref{eq:statbvp}, for $n=2$ this is the Gaussian~curvature.

ii) Proposition~\ref{thm:A}\,iii) with $\mu_0=\mu_1=0$ and $r=\lambda_{1}$ gives $\widehat{\rho}(t,\cdot)\to(A_{1}/\delta)\sin(\pi u)$ in $C^{1}$, and $A_{1}>0$ since $\rho_0>0$ on $(0,1)$ makes the integrand of $A_1$ positive there.
Similar argument as in (1) yields uniform convergence of the generating curves to the curve \eqref{eq:piL} with $B=A_{1}/\delta>0$, and passing to the limit in $|\,\widehat{\rho}\,'_{u}\,|\le1$ gives $B\,\pi\le1$,
i.e., $B\le 1/\pi$, so the limit lies in the extremal family of Proposition~\ref{prop:rbound}, for $B=1/\pi$ the limiting hypersurface closes smoothly at both poles ($|\,\widehat\rho\,'_u\,|=1$ there), and for $B<1/\pi$ it closes with two conical points.
By the above, $\delta\le1$, so \eqref{eq:embedding} applies to $\rho$ (in place of $\,\widehat\rho$\,), and the proof is similar
as for $\delta=1$: we have 
\begin{equation}\label{eq:slopecond}
\max_{u\in[0,1]}\abs{\,\rho'_\infty}
=\max_{u\in[0,1]}\Big|\frac{({\mu_1}/{\delta})\Cn_{r}(u)
-({\mu_0}/{\delta})\Cn_{r}(1-u)}{\Sn_{r}(1)}\Big|\le 1
\end{equation}
for the limiting profile curve, \eqref{eq:limitgen} with $(\mu_0,\mu_1)$, and the slope condition.
\end{proof}

\subsection{Parallel Curved Surfaces}

A \emph{parallel curved (PC) surface} is a surface ${M}^2\subset\Real^{3}$ and a fixed plane $E$ such that at each point of ${M}^2$ at least one principal direction is parallel to $E$, see~\cite{a2}.
PC surfaces carry twisted product metrics, see \eqref{eq:semigeo}, with (positive by the focal condition)
\begin{equation}\label{eq:bclasses}
\rho_{0}(u,v)=1-\rho_{0}(u)k(v),\qquad \text{respectively,}\qquad \rho_{0}(u,v)=1-{a}_{0}(u,v)k(v).
\end{equation}

Let $\gamma:I\to\Real^{2}\subset\Real^{3}$ be a simple closed convex curve, parameterized by arclength,
with curvature $k=k_{\gamma}(v)>0$ and inward unit normal $N(v)$, $r_{\,\mathrm{foc}}:=\min_{\,v}1/k(v)>0$ its focal radius, and $e_3$ the vertical unit vector in $\Real^{3}$.
Let $\rho,h:[0,1]\to\Real$ be smooth arclength
parameterization and $\abs{\rho(u)}<r_{\,\mathrm{foc}}$ for all $u$. Then
 ${\bf r}=\gamma(v)+\rho(u)\,N(v)+h(u)\,e_3$
parameterizes a PC surface, called a \emph{canonical} PC surface, and

\ \
 i)~the induced metric ${\rm d}s^2=\,{\rm d}u^{2}+\big(1-\rho(u)\,k(v)\big)^{2}{\rm d}v^{2}$ belongs to twisted products

\ \
 ii)~$u$-curves are unit-speed geodesics of length $1$

\ \
 iii)~the Gaussian curvature is
$K={\rho''_{u,u}\,k}/({1-\rho\,k})$

\ \
 iv)~with unit normal $\nu=-h'_u N+\rho'_u e_3$, the second fundamental form is diagonal with
$e=h''_{u,u}\,\rho'_u-h'_u\,\rho''_{u,u}$ and $g=-h'_u(u)\,k(v)\big(1-\rho(u)k(v)\big)$, so the coordinate directions are principal, ${\bf r}'_{,v}$ being parallel to the $xy$-plane.

\begin{proposition}[The PRF on PC surfaces] Let $g_{0}$ be induced metric of a canonical PC surface, with $\rho_{0}$ as in \eqref{eq:bclasses} and fixed boundary values $\mu_0(v)=\rho_{0}(0,v)$ and $\mu_1(v)=\rho_{0}(1,v)$. 
Then the {\rm PRF} has a unique global solution, the foliation by  $u$-curves remains geodesic with leaves of length $1$, and for $r<\lambda_{1}$ the metrics converge in every $C^{k}$-norm to~metric
$g_{\infty}=\,{\rm d}u^{2}+\rho_\infty\big(u,r,\mu_0(v), \mu_1(v)\big)^{2}\,{\rm d}v^{2}$
of constant Gaussian curvature $r$, the cases $r\ge\lambda_{1}$ are as in \textup{Proposition~\ref{thm:A}}.
\end{proposition}

\begin{proof} This is standard.
\end{proof}

\begin{theorem}[Realization of {\rm PRF} by PC surfaces] Consider canonical PC surfaces over the fixed base $\gamma$ with evolving profile
\begin{align}\label{eq:pcfamily}
{\bf r}(u,v,t)=\gamma(v)+\rho(t,u)\,N(v)+h(t,u)\,e_3,
\qquad
h(t,u)=\!\!\int_{0}^{u}\!\!\sqrt{1-\rho'_{s}(t,s)^{2}}\,{\rm d}s .
\end{align}

$1)$ {\rm (Rigidity)}. If the induced metrics of a family \eqref{eq:pcfamily} with $\rho'_t\not\equiv0$ satisfy the {\rm PRF}, then either $k={\rm const}$, i.e., $\gamma$ is a circle and the surfaces are rotational, or $r=0$.

$2)$ {\rm (Existence for $r=0$)}. Let $\rho(t,u)$ solve the heat equation $\rho'_{t}=\rho''_{u,u}$ with data $\rho(\cdot,0)={\phi}$, $\rho(t,0)=\alpha_1$, $\rho(t,1)=\alpha_2$, where ${\phi}$ is smooth,
$\abs{\alpha_1},\abs{\alpha_2},\sup\abs{\phi}<r_{\,\mathrm{foc}}$,
and $\sup_{t,u}\abs{\,\rho'_{u}}\le1$ (achievable by the normalization of \textup{Definition~\ref{defn:normalized}}, which scales $\alpha_2$ toward $0$ and so preserves the focal condition). Then

\ \
$2a)$ the focal condition is preserved\textup{:}  by the maximum principle, for all $t$ we have

\begin{equation*}
\sup\nolimits_{\,u}\,\abs{\,\rho(t,u)}\ \le\ \max\big(\sup\abs{\phi},\abs{\alpha_1},\abs{\alpha_2}\big)\ <\ r_{\,\mathrm{foc}}
\end{equation*}

\ $2b)$ a family \eqref{eq:pcfamily} consists of canonical PC surfaces realizing the {\rm PRF} with~$r=0$.

\ \ $2c)$ the limit is a flat canonical PC surface with
profile $\rho_{\infty}(u)
=\alpha_1+(\alpha_2-\alpha_1)\,u$.
\end{theorem}

\begin{proof}
1) The induced metric of \eqref{eq:pcfamily} is
${\rm d}u^{2}+b^{2}\,{\rm d}v^{2}$ with $b=1-\rho(t,u)k(v)$.
The~PRF is equivalent to $\rho'_{t}=\rho''_{u,u}+r\,b$, thus, (using $k(v)>0$)
\begin{equation}\label{eq:pcansatz2}
\rho'_{t}-\rho''_{u,u}-r\,b=-{r}/{k(v)} .
\end{equation}
The LHS is independent of $v$, hence so is the RHS, which forces $r=0$ or $k={\rm const}$.

2) For $r=0$, \eqref{eq:pcansatz2} reduces to $\rho'_{t}=\rho''_{u,u}$, so the induced metrics satisfy the PRF equation with $r=0$, this proves the realization claim in 2b) once the surfaces are shown to be well defined.
2a) is the maximum principle for the heat equation: both $\max$ and $\min$ of $\rho$ over $[0,1]\times[0,T]$ are attained on the parabolic boundary, where $\abs{\rho}\le
\max(\sup\abs{{\phi}},\abs{\alpha_1},\abs{\alpha_2})<r_{\mathrm{foc}}$, note that for $r=0$ there is no zeroth-order term and the principle applies directly. Together with $\sup\abs{\,\rho'_{u}}\le1$, the map \eqref{eq:pcfamily} is well defined, and $1-b\,k\ge1-\sup\abs{b}\,\max k>0$, so it is an immersion with the stated first fundamental form.
For 2c), Theorem~\ref{T-stabnonstat} with $r=0<\lambda_{1}$ gives $\rho(t,\cdot)\to\rho_\infty(\cdot\,,0,\alpha_1,\alpha_2)$ in $C^{1}$, and the height functions converge uniformly,
the limit profile is linear, so the limiting surface has $K=k\,(\rho_{\infty})''_{u,u}/(1-k\,\rho_{\infty})\equiv0$.
\end{proof}

\section{The PRF on Weak Almost Contact Metric Manifolds}
\label{sec:04}

Contact Riemannian geometry is a significant field in both mathematics and physics. In 2022, Rovenski and Wolak introduced metric structures on a smooth manifold that generalize the almost contact, $K$-contact, Sasakian, cosymplectic, etc. metric structures.
Such so-called ``weak" structures allow us to take a fresh look at the theory of classical structures and find new applications, see \cite{RWo-2,rov-117,rov-126}.

A Riemannian manifold $(M^{2n}, g)$ equipped with a non-degenerate skew-symmetric (1,1)-tensor $f$
is called a \textit{weak Hermitian manifold}.
This manifold is called \textit{weak K\"{a}hlerian} if ${\nabla} f=0$, where $\nabla$ is the Levi-Civita connection.

A \textit{weak almost contact metric structure} on a
manifold $M^{2n+1}$ is a set $({f},Q,\xi$, $\eta,g)$, where $g$ is a Riemannian metric on $M$, ${f}$ is a skew-symmetric (1,1)-tensor of rank $2n$, $\xi$ is a unit (Reeb) vector field, $\eta$ is a 1-form such that ${\eta}({\xi})=1$, satisfying

\begin{equation*}
g({{f}} X,{{f}} Y)= g(X,Q\,Y) -{\eta}(X)\,{\eta}(Y)
\end{equation*}

and $Q:=-{{f}}^2 + {\eta}\otimes {\xi}$ is a self-adjoint (1,1)-tensor field (for almost contact metric manifolds, $Q=I$ is true).
The~self-adjoint (1,1)-tensor $\widetilde Q:=Q-I$
measures the difference between weak and classical metric structures.
A~\textit{weak contact metric structure} is a weak almost contact metric structure satisfying the condition
${\rm d}\eta = F$, where $F(X,Y)=g(X,fY)$,
in this case, the Reeb vector field $\xi$ is geodesic.

\subsection{Weak \texorpdfstring{$K$}{K}-Contact Manifolds}

A weak contact metric structure is called a \textit{weak $K$-contact structure} if
${\xi}$ is Killing (an infinitesimal isometry), i.e., the Lie derivative of metric vanishes: ${\cal L}_{\xi}\,g=0$.
We generalize the fact that a Riemannian manifold with a unit Killing vector field $\xi$ such that $R_{\xi}X={X}$ for ${X}\in\ker\eta$ is a $K$-contact manifold.

\begin{theorem}[\cite{rov-117}]\label{prop2.1b}
A Riemannian manifold $(M^{\,2n+1},g)$ admitting a unit Killing vector field $\xi$ such that the $(1,1)$-tensor $\nabla\xi$ has constant rank $2n$ and $\xi$-sectional curvature is positive is a weak $K$-contact manifold $M(f,Q,\xi,\eta,g)$ with $\eta = g(\cdot\,, \xi)$
and the following structural tensors: $f=-\nabla\,\xi$
and $Q {X} = R_{\xi}X\,$ for ${X}\in\ker\eta$.
\end{theorem}

The PRF prescribes $K_{\rm mix}$ of a foliation
and has recently been used to obtain dynamical characterizations of $K$-contact structures,
see \cite{rov-117}.
For $K$-contact structures, $\Ric^\perp_{g}=1$vand the metric is a fixed point of the normalized PRF with $r=1$.

\begin{theorem}[\cite{rov-117}]
\label{T-PRF-g2b}
Let $M^{2n+1}(f_0,Q_{\,0},\xi,\eta,g_0)$ be a weak $K$-contact manifold.
Then, there exist ${{\cal D}}$-adapted weak $K$-contact structures $(f_t,Q_{\,t},\xi,\eta,g_t)$ on $M$
such that $g_t$ -- a solution of the {\rm PRF} -- converges exponentially fast, as $t\to-\infty$, to a limit metric
that provides a $K$-contact structure.
\end{theorem}

\begin{example}\label{Ex-K-S}
\rm
Let $M^{2n+1}(f_0,Q_0,\xi,\eta,g_0)$ be a weak $K$-contact manifold
with $Q|_{\,\cal D} =\lambda\,I_{\,\cal D}$ for some differentiable function
$0<\lambda\not\equiv1$ on~$M$.
From Theorem~\ref{T-PRF-g2b}, applying the theorem of uniqueness of a solution to ODE
\begin{equation}\label{E-PRF-fK}
 \partial_t R_\xi = 4\,R_\xi\,(R_\xi -\,{I}_{\cal D})
\end{equation}
we have a $\cal D$-conformal solution $R_\xi(t)
=\mu(t)\,I_{\cal D}$.
Since $\partial_t I_{\cal D}=0$, we have also that  ODE
${\rm d}\mu/{\rm d}t=4\mu\,(\mu-1)$ with $\mu(0)=\lambda>0$.
Thus, $\mu(t)=\dfrac{1}{1-C\,\exp(4t)}\ (t\in\mathbb{R})$, where $C=1-\lambda^{-1}\ne0$.
Since our PRF reduces to ODE, by the uniqueness of a solution, we~have $\partial_t g\,\!^\bot_t = 0$ and
$g_t=g^\bot_0+\rho(t)\,g^\top_0$, where
$\rho(0)=1$ and
\begin{eqnarray*}
 &\partial_t g\,\!^\top_t \!= 2(1{-}\mu(t))g^\top_t
 {\implies} 
 \dfrac{\rm d}{{\rm d}t}\log\rho(t)
 \!= {-}\dfrac{2\,C\exp(4t)}{1{-}C\exp(4t)}
 {\implies}
 \rho(t) \!= \Big(\dfrac{1{-}C\exp(4t)} {1-C}\Big)^{1/2}.
\end{eqnarray*}
Therefore, $f_t=f_0/\rho(t)$ and $Q_{\,t}|_{\,\mathcal D} =\dfrac{1}{1-C\exp(4t)}\,I_{\mathcal{D}}$.
Note that $\lim\limits_{t\to-\infty}Q_t=I_{\cal D}$.
By~Theorem~\ref{T-PRF-g2b},
$(f_t,Q_{\,t},\xi,\eta,g_t)$ $(t\in\mathbb{R})$ is a smooth family of weak $K$-contact structures on $M$
-- a solution of the PRF
-- that converges exponentially fast, as $t\to-\infty$, to a $K$-contact structure
$(f_0/{\sqrt\lambda}, \xi,\eta, g^\bot_0+\sqrt{\lambda}\,g^\top_0)$.
\end{example}

We obtain an extrinsic realization of the PRF on an ellipsoid  of $\mathbb{R}^4$
with the weak $K$-contact structure, its limit is a unit three-sphere with the $K$-contact structure.

\begin{example}
Let $M\subset\mathbb{R}^{2n+2}$ with $n\ge1$ be an ellipsoid with induced metric $g$
\[
 M = \Big\{u=(u_1,\ldots,u_{2n+2})\in\mathbb{R}^{2n+2}; \qquad \sum\nolimits_{\,i=1}^{n+1} u_i^2 + a\sum\nolimits_{\,i=n+2}^{2n+1} u_i^2 = 1\Big\}
\]
where $0<a={\rm const}\ne1$. The sectional curvature of $(M,g)$ is positive, and
\[
 \xi(u) = (-u_2, u_1, \ldots , -u_{n+1}, u_{n}, -\sqrt a\,u_{n+3}, \sqrt a\,u_{n+2}, \ldots , -\sqrt a\,u_{2n+2}, \sqrt a\,u_{2n+1})
\]
is a Killing vector field on $\mathbb{R}^{2n+2}$.
(Indeed, there is a constant block-diagonal with $2\times2$ skew-symmetric blocks matrix $J$ such that $\xi(u)=Ju$.
The integral curves of $\xi(u)$ solve
$({\rm d}/{\rm d}t)u=Ju$, hence, $u(t)=\exp(tJ) u(0)$, where $\exp(tJ)$ is an orthogonal matrix.
Thus, the flow $\Psi_t(u)=\exp(tJ) u$ consists of rotations, so each $\Psi_t$ is an isometry).
Note that $\xi$ is tangent~to~$M$,
and $\xi$ is a unit Killing vector field on $(M,g)$, which defines a weak $K$-contact structure, see
Theorem~\ref{prop2.1b}.
\end{example}

The following result (for $n=1$) can easily be extended for $n>1$.

\begin{theorem}
Let $(f,Q,\xi,\eta,g)$ be a weak $K$-contact structure on an ellipsoid
$S^3=\big\{u_1^2 + u_2^2 + a_0(u_3^2 + u_4^2) = 1\big\}\subset\mathbb{R}^{4}$ with a unit Killing vector field
$\xi = (-u_2,u_1,-\!\sqrt a_0\,u_{4}$, $\sqrt a_0\,u_{3})$ and $Q|_{\,\ker\eta}=\lambda\,I_{\,\ker\eta}$ for a smooth function $\lambda>0$ on~$S^3$.
Then we get an extrinsic realization of the {\rm PRF}
on elipsoids $S^3_t=\big\{u_1^2+u_2^2+a_t(u_3^2+u_4^2)=1 \big\}\subset\mathbb{R}^{4}$ with $a_t=1-(1-a_0)\exp(4\,t)$ and $t\ge0$, and the limit, as $t\to-\infty$, is a unit sphere $\mathbb{S}^3\subset\mathbb{R}^{4}$ equipped with a $K$-contact structure $(\hat f,\hat \xi,\hat \eta,\hat g)$ and $\hat \xi = (-u_2, u_1, -u_{4}, u_{3})$.
\end{theorem}

\begin{proof}
Using Hopf parametrization
$u_1=\cos\theta\cos\phi$,
$u_2=\cos\theta\sin\phi$,
$u_3=a_0^{-1/2}\sin\theta\cos\psi$,
$u_4=a_0^{-1/2}\sin\theta\sin\psi$)
of this ellipsoid,
we find its metric tensor

\begin{equation*}
{\rm d}s^2= \big(\sin^2\theta+a_0^{-1}\cos^2\theta\big)\,{\rm d}\theta^2
 +\cos^2\theta\,{\rm d}\phi^2
 +a_0^{-1}\sin^2\theta\,{\rm d}\psi^2.
\end{equation*}

For any tangent vector $X$ such that
$g(X,\xi)=0$, the Gauss equation, and \eqref{E-PRF-fK},
give
 $R_{\,\xi}X=g(A\xi,\xi)AX -g(AX,\xi)A\xi$,
where $A$ is the Weingarten operator.
In fact, for this ellipsoid of revolution,
$X$ and $\xi$ are principal directions.
By the above,
 $R_{\,\xi}X = \|u\|^{-2}\,X
 =(1+(1-a_0)(u_3^2+u_4^2))^{-1}X$,
where the factor varies with both
$a_0$ and the point $(u_i)$ on the ellipsoid.
Using $QX=R_{\,\xi}X$ (Theorem~\ref{prop2.1b}) where $X\in\ker\eta$
and Example~\ref{Ex-K-S}, we get
 $1+(1-a_t)(u_3^2+u_4^2)= 1-(1-\lambda^{-1})\,\exp(4\,t)$.
For $a_t=1-(1-a_0)\exp(4\,t)$ and $t\in\mathbb{R}$,
we obtain
$\lambda=(u_3^2+u_4^2+1)^{-1}$
 and
 $QX=(1+\exp(4\,t)(u_3^2+u_4^2))^{-1}X\
 (X\in\ker\eta)$.
We get an extrinsic realization of the PRF on the weak $K$-contact ellipsoid. By~Theorem~\ref{T-PRF-g2b}, the limit, as $t\to-\infty$, is a unit sphere $\mathbb{S}^3$ with a $K$-contact~structure.
\end{proof}

\subsection{Weak \texorpdfstring{$\beta$}{beta}-Kenmotsu Manifolds}

In 1972, Kenmotsu introduced and studied a special class of almost contact metric manifolds,
which is characterized in terms of warped products
of the real line and a K\"{a}hler manifold
with the warping function $\rho=C\exp(u)$.
In 1991, Olszak characterized in terms of warped products (with $\rho=C\exp(\beta\,u)$) a class of almost contact metric manifolds, known as $\beta$-Kenmotsu manifolds for $\beta={\rm const}>0$, given~by
\begin{equation}\label{2.3}
 (\nabla_{X}\,f)Y=\beta\,\big(\,g(fX, Y)\,\xi-\eta(Y)\,fX\,\big).
\end{equation}
\textit{Weak $\beta$-Kenmotsu manifolds} were introduced, see \cite{rov-126}, as weak almost contact metric manifolds satisfying \eqref{2.3}, where $\beta\not\equiv0$ is a smooth function on $M$, and if $n=1$, then  ${\rm d}\beta\wedge\eta=0$ is required.
If $\beta\equiv1$, then $M^{2n+1}(f,Q,\xi,\eta,g)$ is called a \textit{weak Kenmotsu manifold},
and for $Q=I$ we get a \textit{Kenmotsu manifold}.
By \eqref{2.3},
$\nabla_X\,\xi=\beta\,\big(\,X-\eta(X)\xi\,\big)$
thus, for weak $\beta$-Kenmotsu manifolds,
$\xi$ is geodesic~($\nabla_\xi\,\xi=0$).

\begin{theorem}[see \cite{rov-126}]\label{T-2.1}
Every point of a weak $\beta$-Kenmotsu manifold $M({f},Q, \xi,\eta,g)$ has a neighborhood $U$ that is isometric to a twisted product $(-\varepsilon,\varepsilon)\times_\rho\bar M$
$($a warped product if $\nabla\beta$ is parallel to $\xi$
with the leaves tangent to $\xi$, where $\bar M$
has a weak Hermitian (weak K\"{a}hler, respectively) structure and the following is true
\begin{equation}\label{E-f-warp}
 \rho'_u/\rho=\beta.
\end{equation}
For a simply connected $M$ and geodesically complete $\xi$, the isometry is~global.
\end{theorem}

Consider a (weak) $\beta$-Kenmotsu structure on a warped product $[0,1]\times_\rho\mathbb{S}^2(1)$ with the (weak) K\"{a}hler structure on $\mathbb{S}^2(1)$, and realize it on a rotational hypersurface in $\mathbb{R}^4$ with profile
curve
$(\rho(u),h(u))\ (u\in[0,1])$, where $(\rho'_u)^2+(h'_u)^2=1$.
By~Theorem~\ref{T-2.1}, $\nabla\beta$ is parallel to the Reeb field $\xi$, i.e., $\beta$ is a function of~$u$.
By \eqref{E-f-warp}, we have
\begin{align*}
\rho(u)&=\rho_0\exp\big(\int_0^u\beta(s)\,{\rm d}s\big),\\
h(u)&=\int_0^u\sqrt{1-\beta^2(s)\,\rho_0^2
\exp\,\Big(2\int_0^s\beta(\tau)\,{\rm d}\tau\Big)}\,{\rm d}s.
\end{align*}

For~$\beta\equiv1$, the profile curve is
 $(\rho(u),h(u))= \Big(  \rho_0 \exp\,(u),\,
\int_0^u \sqrt{1-\rho_0^2 \exp\,(2s)}\,{\rm d}s\Big)$.

Let the family of warped products $M^3_t=[0,1]\times_{\rho(t,u)}\mathbb{S}^2$, $t\ge0$,
with $\rho(0,u)$ given in i), see Theorem~\ref{T-stabnonstat},
carries a (weak) $\beta$-Kenmotsu structure with $\xi=\partial_u$ and $\beta$ depending on $(t,u)$, and $g(t)$ modified by the PRF \eqref{E-flow-a}.
Let us realize $M^3_t$ as hypersurfaces of revolution $S^3_t$ in $\mathbb{R}^4$ with profile curves $(\rho(t,u),h(t,u))$, where
the height function $h(t,u)$ is determined (up to an additive constant in $u$) by the arc-length condition,
i.e., $h'_u(t,u) = \sqrt{1 - (\rho'_u)^2(t,u)}$.
The rotational hypersurface~closes smoothly at both ends
$u=0$ and $u=1$ (i.e., is diffeomorphic to $\mathbb{S}^3$)
if and only if
\begin{eqnarray*}
\rho(t,0) = 0,\ \ \rho(t,1) = 0,\quad
\rho'_u(t,0) = 1,\ \ \rho'_u(t,1) = -1, \quad
 h'_u(t,0) = 0,\ \ h'_u(t,1) = 0.
\end{eqnarray*}
By \eqref{E-f-warp}, $\beta=\rho'_u/\rho$ is a function of $(t,u)$ that defines a $\beta$-Kenmotsu structure on~$M^3_t$.

The PRF modifies $g(t)$, as $t\to\infty$, to a limit metric
$g_\infty = {\rm d}u^2 + \rho_\infty^2(u)\,g_{\mathbb{S}^2}$
with $K_{\rm mix}\equiv r$. The function $\rho_\infty(u)$
satisfies the ODE and is given by \eqref{E-stat-rho}.
The profile curve of the limiting rotational hypersurface
$S^3_\infty$ in $\mathbb{R}^4$ is given by
$\bigl(\rho_\infty(u),h_\infty(u)\bigr)$,
where
$h_\infty(u) = \int_0^u
\sqrt{1 - (\rho_\infty')^2(s)}\,{\rm d}s$.
The limiting (weak) $\beta$-Kenmotsu structure
has $\beta_\infty=(\rho_\infty)'_u/\rho_\infty$.
For $r=\pi^2$, we get $\rho_\infty = C\sin(\pi u)$, see Example~\ref{Ex-tilder}, then using \eqref{E-f-warp} we find $\beta_\infty(u)=\pi\cot(\pi u)$.

\section{Conclusion}

We studied extrinsic representations of the partial Ricci flow
(PRF) on manifolds with geodesic foliations. The constructions encompass
compact rotational hypersurfaces in Euclidean space, including weak $\beta$-Kenmotsu examples,
parallel curved surfaces, and ellipsoidal models
carrying weak $K$-contact structures.

The study of extrinsic representations of the PRF presented in this paper can be extended for totally geodesic foliations, e.g., twisted products and weak metric $f$-contact manifolds (as well as to spherical foliations such as in \cite{GT-2026}) with applications to the calculus of moving surfaces and mathematical physics.

\section*{Acknowledgments}
The authors thank Professor D.\,Blanc for helpful comments on the paper.

\end{document}